\documentclass{article}
\usepackage{amsmath}

\newtheorem{theorem}{Theorem}[section]

\newtheorem{lemma}[theorem]{Lemma}
\newtheorem{cor}[theorem]{Corollary}
\newtheorem{proposition}[theorem]{Proposition}
\newtheorem{definition}[theorem]{Definition}
\newtheorem{example}[theorem]{Example}

\newtheorem{remark}[theorem]{Remark}
\newcommand{\be}{\begin{equation}}
\newcommand{\ee}{\end{equation}}
\newcommand{\ben}{\begin{enumerate}}
\newcommand{\een}{\end{enumerate}}

\def\beq{\begin{equation}}
\def\eeq{\end{equation}}

\newcommand{\qed}{\hspace*{\fill}Q.E.D.}  

\title{{\bf\LARGE On the Berwald-Weyl Curvature }}
\author{Zhongmin Shen\footnote{supported in part by NSF China (NSFC No. 11671352) }\;
and Liling Sun\footnote{corresponding author, supported by NSF China (NSFC No. 12101259), Natural Science Foundation of Fujian Province(No. 2021J05173), Foundation of Fujian Educational Committee(No. JAT200274) and Research Foundation of Jimei University(No. ZQ2022032)}}
\date{}

\begin{document}

\maketitle
\begin{abstract}
In this paper, we study the Berwald-Weyl curvature which is defined for a spray/Finsler metric with a volume form. We obtain some expressions for the Berwald-Weyl curvature. This quantity is a projective invariant with respect to a fixed volume form.
We prove that for any spray of scalar curvature on a  manifold of dimension $n\geq 3$,  the  Berwald-Weyl curvature  vanishes with respect to any volume form.  We also show that for any Finsler metric of constant Ricci curvature and constant S-curvature, the Berwald-Weyl curvature vanishes with respect to the Busemann-Hausdorff volume form. This study leads to a new notion of  BWeyl-flat sprays/Finsler metrics.


\end{abstract}

\section{Introduction}
In Finsler geometry, there are several well known  projective invariants  such as Douglas curvature and Weyl curvature.  Finsler metrics with vanishing  Weyl curvature are called Weyl metrics, that are
 exactly Finsler metrics of scalar flag curvature (\cite{Z.Sza}). In dimension $n\geq 3$, a spray/Finsler metric is locally projectively flat if and only if the Weyl curvature and the Douglas curvature vanish (\cite{Dou}) . In dimension two, the Weyl curvature always vanishes. There is another projective invaruant, the Berwald-Weyl curvature. It is known that a two-dimensional spray/Finsler is locally projectively flat if and only if the Berwald-Weyl curvature and the Douglas curvature vanish (\cite{L.Be}). In this paper, we shall discuss the Berwald-Weyl curvature and its generalzation to higher dimensions.

A spray on a manifold $M$ is a vector field  on the tangent bundle  $TM$ in the form $\mathbf{G}=y^{i}\frac{\partial}{\partial x^{i}}-2G^{i}\frac{\partial}{\partial y^{i}} $. Fix a local coordinate system, one can construct a local projective spray $\bar{\bf G}$ by setting $\bar{G}^i =G^{i}-\frac{1}{n+1}\frac{\partial  G^m}{\partial y^m}y^{i}.$ Define a tensor ${\bf W}^o$ by
\[ W^o_k := \bar{R}_{\parallel k}-\frac{1}{2}\bar{R}_{\cdot k\parallel m}y^{m},
\]
where the covariant derivative $``\parallel"$ is taken respect to $\bar{\bf G}$ (\cite{S1}) and $\bar{R}=\frac{1}{n-1}{\bf Ric}_{\bar{G}}$ is the Ricci scalar of $\bar{\bf G}$. It is proved that in dimension two, ${\bf W}^o$ is a globally defined quantity (\cite{L.Be} \cite{S1}).
In higher dimensions, the
above expression does not define a global quantity.
To avoid this problem, one can use  the globally defined associated projective spray $\hat{\bf G}$ using a volume form \cite{S1}.
Given  a volume form $dV$ on $M$, we have the S-curvature ${\bf S} = {\bf S}_{(G, dV)}$ (\ref{621}). With the S-curvature, we can construct a new spray by
\begin{equation}
 \hat{\mathbf{G}}:=\mathbf{G}+\frac{2\mathbf{S}}{n+1}\mathbf{Y},\label{hatG}
\end{equation}
where   $\mathbf{Y}=y^i\frac{\partial}{\partial y^i}$ denotes the canonical vertical field on $TM$. This spray is projectively invariant  with respect to  a fixed volume form.
The spray $\hat{\bf G}$ is called the projective spray of $({\bf G}, dV)$. Let $\hat{R}=\frac{1}{n-1}{\bf Ric}_{\hat{G}}$ is the Ricci scalar of $\hat{\bf G}$. The Berwald-Weyl curvature ${\bf W}^o = W^o_k dx^k$ of $({\bf G}, dV)$ is defined by
\begin{equation}\label{7261}
W^{o}_{k}=\hat{R}_{\parallel k}-\frac{1}{2}\hat{R}_{\cdot k\parallel m}y^{m},
\end{equation}
where the covariant derivative $``\parallel"$ is taken respect to $\hat{\mathbf{G}}$ (\cite{S1}). The Berwald-Weyl curvature is projectively invariant with respect to a fixed volume form. \\

 In general,  the spray $\hat{\mathbf{G}}$ defined in (\ref{hatG})  depends on the volume form $dV$ and so does the Berwald-Weyl curvature ${\bf W}^o$.
In $\dim M\geq3$, ${\bf W}^o$  can be expressed in terms of the Weyl curvature.
\begin{proposition}\label{W0kN3}
Given a spray space $(\mathbf{G},dV)$ on the manifold $M$ of dimension $n\geq 3$, the Berwald-Weyl curvature $\mathbf{W}^{o}$ can be written as
\begin{equation}\label{d9d}
W^{o}_{k}=\frac{1}{n-2}W^{m}_{~k||m},
\end{equation}
where the covariant derivative $``\parallel"$ is taken respect to $\hat{\mathbf{G}}$ and $W^{i}_{k}$ is the Weyl curvature given in \eqref{d7}.
\end{proposition}

It is well known that a spray is of scalar curvature if and only if its Weyl curvature $\mathbf{W}$ vanishes, that is $W^i_{~k}=0$.
By this fact and \eqref{d9d}, we conclude immediately the following
\begin{theorem}\label{e88}
If $\mathbf{G}$ is a spray of scalar curvature on a manifold $M$ of $\dim  M\geq 3$,
then its Berwald-Weyl curvature $\mathbf{W}^{o}=0$ for any given volume form $dV$.
\end{theorem}

Since the Berwald -Weyl curvature is defined via the S-curvature, it is natural to express ${\bf W}^o$ in terms of the geometric quantities of ${\bf G}$ and the S-curvature with respect to the Berwald connection of ${\bf G}$.

\begin{proposition}\label{BWC1}
Given a spray space $(\mathbf{G},dV)$ on the manifold $M$ of dimension $n$, the Berwald-Weyl curvature $\mathbf{W}^{0}$ can be expressed as
\begin{equation}\label{d99}
W^{o}_{k}=R_{|k}-\frac{1}{2}R_{\cdot k|m}y^{m}-\frac{1}{n+1}\chi_{k|m}y^{m}-\frac{1}{n+1}W^{m}_{~k}\mathbf{S}_{\cdot m},
\end{equation}
where  $R=\frac{1}{n-1}{\bf Ric}_G$ is the Ricci scalar of ${\bf G}$,    $\chi_k$ is the non-Riemannian quantity given in (\ref{chi}),
 $W^{i}_{k}$ is the Weyl curvature given in \eqref{d7},  ${\bf S} = {\bf S}_{(G, dV)}$ is the S-curvature of $({\bf G}, dV)$ given in (\ref{621}) and the horizontal covariant derivatives $``|"$ are taken with respect to $\mathbf{G}$.
\end{proposition}

We know that if a Finsler metric is of constant flag curvature $k$, then $R =k F^2$, $\chi_k =0$ and $W^i_k =0$.
From (\ref{d99}), one can see that  $W^{o}_{k}=0$ for any volume form.

In dimension $n=2$, since $\bf{W}=0$. By this fact and \eqref{d99}, we see that ${\bf W}^o$ is independent of the volume form. This fact is known  \cite{S1}.

\begin{cor}\label{cor1.4} For a two-dimensional spray $({\bf G}, dV)$, the Berwald-Weyl curvature ${\bf W}^o$ is independent of $dV$, given by
\begin{equation}\label{d99*}
W^{o}_{k}=R_{|k}-\frac{1}{2}R_{\cdot k|m}y^{m}-\frac{1}{3}\chi_{k|m}y^{m}.
\end{equation}
Hence it is a projective invariant of ${\bf G}$.
\end{cor}

Every Finsler metric  $F$ induces a spray  ${\bf G}_F$. Thus results on sprays are still true.
However, for a Finsler metric $F$,  we sometimes restrict the volume form $dV$  to its Busemann-volume form $dV_{BH}$. In this case, the S-curvature of $F$ is the S-curvature of $({\bf G}_F, dV_{BH})$. A Finsler metric
 $F$ is of isotropic S-curvature if the S-curvature
\[\mathbf{S}=(n+1)cF,\]
where $c=c(x)$ is a scalar function on $M$.
$F$ is of constant S-curvature if the $c$ is a constant. $F$ is called an Einstein metric if
\begin{equation*}
\mathbf{Ric}=(n-1)\sigma F^{2},
\end{equation*}
where $\sigma=\sigma(x)$ is a scalar function on $M$. In particular, a Finsler metric $F$ is of constant Ricci curvature if $\sigma$ is a constant.\\

For a Finsler manifold $(M, F)$, if the S-curvature is constant with respect to  the Busemann-Hausdorff volume form $dV=dV_{BH}$, then $\chi_k=0$.  Further, if the Ricci curvature is constant, then by  \eqref{d99}, the Berwald-Weyl curvature vanishes

\begin{theorem}\label{E-RC}
Let $(M, F)$ be a Finsler manifold. If $F$ is of constant Ricci curvature and constant S-curvature with respect to $dV=dV_{BH}$, then ${\bf W}^o=0$ .
\end{theorem}

It is known that every Einstein Randers metric on a manifold $M$   is of constant S-curvature\cite{BR}. Hence, we obtain the following
\begin{cor}\label{R-RC}
For any Randers metric of constant Ricci curvature, ${\bf W}^{o}=0$ for the $dV=dV_{BH}$.
\end{cor}

We know that there exist Randers metrics of constant S-curvature and constant Ricci curvature which are not of scalar curvature \cite{CS2}. Thus there exist Finsler metrics with vanishing Berwald-Weyl curvature which are not of scalar curvature. Below is another non-trivial example.

\begin{example}
Let $\alpha_{1}=\sqrt{a_{ij}y^{i}y^{j}}$ and $\alpha_{2}=\sqrt{\bar{a}_{ij}y^{i}y^{j}}$ be two Ricci-flat Riemannian
metrics on the manifolds $M_{1}$ and $M_{2}$, respectively. Consider the following 4-th root metric
\begin{equation*}
F:=\sqrt[4]{\alpha_{1}^{4}+2c\alpha_{1}^{2}\alpha_{2}^{2}+\alpha_{2}^{4}},
\end{equation*}
where $0 < c \leq 1$ is a constant.  It is a Riemannian metric when $c =1$.
This is a Ricci-flat   Berwald metric on $M:=M_{1}\times M_{2}$.  Note that   the S-curvature ${\bf S}_{(F, dV)}=0$ for the Busemann-Hausdorff volume form $dV=dV_{BH}$. Thus ${\bf W}^o=0$.
\end{example}

\section{Preliminaries}
Let $M$ be an $n$-dimensional manifold and $TM$ its tangent bundle. We denote by $TM_{0}=TM\setminus\{0\}$ the slit tangent bundle. Local coordinates on the base manifold $M$ will be denoted by $(x^{i})$, while the  standard local coordinates on $TM$ or $TM_{0}$ will be denoted by $(x^{i} ,y^{i})$.

A spray $\mathbf{G}$ on $M$ is a smooth vector field on $TM_{0}$ expressed in a standard local coordinate system $(x^{i},y^{i})$ in $TM$ as follows
\begin{equation}
\mathbf{G}=y^{i}\frac{\partial}{\partial x^{i}}-2G^{i}\frac{\partial}{\partial y^{i}},
\end{equation}
where $G^i=G^{i}(x, y)$ are the local functions on $TM$ satisfying
$
G^{i}(x, \lambda y)=\lambda^{2}G^{i}(x, y) ~\text{for}~\forall\lambda>0.
$
Every Finsler metric $F=F(x,y)$ on a manifold induces a spray with the \emph{geodesic coefficients}
\begin{equation}\label{631}
G^{i}=\frac{1}{4}g^{il}\{[F^{2}]_{x^{k}y^{l}}y^{k}-[F^{2}]_{x^{l}}\},
\end{equation}
where $g^{ij}=(g_{ij})^{-1}$.

The vertical tangent bundle of $TM_0$ is expressed by
\[\mathcal{V}TM:=\text{span}\{\frac{\partial}{\partial y^{1}},\cdots \frac{\partial}{\partial y^{n}}\}.\]
Put
\begin{equation}\label{Hframe}
\frac{\delta}{\delta x^{k}}:=\frac{\partial}{\delta x^{k}}-N^{l}_{k}\frac{\partial}{\delta y^{l}},
\end{equation}
where $N^{i}_{j}=\frac{\partial G^{i}}{\partial y^{j}}.$
Then the horizontal tangent bundle of $TM_0$ is expressed by
\begin{equation*}
\mathcal{H}TM:=\text{span}\{\frac{\delta}{\delta x^{1}},\cdots \frac{\delta}{\delta x^{n}}\}.
\end{equation*}
The tangent bundle $T(TM_{0})$ of $TM_{0}$ is splitted into the direct sum $T(TM_{0})=\mathcal{V}TM\oplus\mathcal{H}TM$.

Set $\omega^{i}=dx^{i} $ and $\omega^{n+i}:=dy^{i}+N^{i}_{j}dx^{j}$.
The local connection 1-forms of the Berwald connection are given by
\begin{equation*}
\omega^{\ i}_{j}:=\varGamma^{i}_{jk}dx^{k},
\end{equation*}
where
\begin{equation*}
\varGamma^{i}_{jk}=\frac{\partial N^{i}_{j}}{\partial y^{k}}=\frac{\partial^{2}G^{i}}{\partial y^{j}\partial y^{k}}.
\end{equation*}
We have
\begin{equation*}
d\omega^{i}=\omega^{j}\wedge\omega^{\ i}_{j}.
\end{equation*}
Then we obtain the curvature form:
\begin{equation}
\Omega^{\ i}_{j}=d\omega^{\ i}_{j}-\omega^{\ l}_{j}\wedge\omega^{\ i}_{l}.
\end{equation}
Put
\begin{equation}\label{d2}
\Omega^{\ i}_{j}=\frac{1}{2}R^{~i}_{j~kl}\omega^{k}\wedge \omega^{l}-B^{~i}_{j~kl}\omega\wedge \omega^{n+l}.
\end{equation}
The two curvature tensors $R^{~i}_{j~kl}$ and $B^{~i}_{j~kl}$ are called \emph{Riemann curvature tensor} and \emph{Berwald curvature tensor}, respectively.

The Riemann curvature $\mathbf{R}_{y}=R^{i}_{~k}dx^{k}\otimes \frac{\partial}{\partial x^{i}}|_{x}:T_{x}M \rightarrow T_{x}M$ of a spray
$\bf{G}$ is a family of linear maps on tangent spaces which is defined by
\begin{equation}\label{611}
R^{i}_{~k}=y^{j}R^{~i}_{j~kl}y^{l}.
\end{equation}
In a standard local coordinates, it is expressed by
\begin{equation*}
R^{i}_{~k}=2\frac{\partial G^{i}}{\partial x^{k}}-y^{j}\frac{\partial^{2} G^{i}}{\partial x^{j}\partial y^{k}}+2G^{j}\frac{\partial^{2} G^{i}}{\partial y^{j}\partial y^{k}}-\frac{\partial G^{i}}{\partial y^{s}}\frac{\partial G^{s}}{\partial y^{k}}.
\end{equation*}
They are related by
\begin{equation}\label{598}
R^{~i}_{j~kl}=\frac{1}{3}(R^{i}_{~k\cdot l}-R^{i}_{~l\cdot k})_{\cdot j}
\end{equation}
Here and hereafter,  $``\cdot"$ denotes the vertical derivatives with respect to $y$.
For instance, $f_{\cdot k}=\frac{\partial f}{\partial y^{k}},~~f_{\cdot k\cdot l}=\frac{\partial^{2}f}{\partial y^{k}\partial y^{l}}$, etc.
By \eqref{611} and \eqref{598}, we see that the  curvature tensors $R^{i}_{~k}$ and $R^{~i}_{j~kl}$ can represent each other.
For this reason, they are all called \emph{Riemann curvature tensor} if there is no confusion.
Put
\begin{equation}
R^{~i}_{~kl}:=y^{j}R^{~i}_{j~kl}.
\end{equation}
By  \eqref{598}, one has
\begin{equation}\label{R3eq}
R^{~i}_{~kl}=\frac{1}{3}(R^{i}_{~k\cdot l}-R^{i}_{~l\cdot k}).
\end{equation}
The Ricci curvature is the trace of $\mathbf{R}_{y}$ defined by
\begin{equation}\label{692}
{\bf Ric}:=R^{m}_{~m}=y^{j}R^{~m}_{j~ml}y^{l}
\end{equation}
and  $R:=\frac{1}{n-1}\bf Ric$ is called the Ricci scalar of $\mathbf{G}$.

Define a linear map $\mathbf{T}_{y}:T_{x}M\rightarrow T_{x}M$ by
\[
\mathbf{T}_{y}(u)=T^{i}_{~j}(x,y)u^{k}\frac{\partial}{\partial x^{i}},~~u=u^{k}\frac{\partial}{\partial x^{k}},\]
where
\begin{equation}\label{Tij}
T^{i}_{j}:=R^{i}_{~j}-(R\delta^{i}_{~j}-\frac{1}{2}R_{\cdot j}y^{i}).
\end{equation}
A spray is said to be of {\it isotropic curvature} if ${\bf T}=0$.

Given a volume form  $dV=\sigma(x) dx^{1}\cdots dx^{n}$ on a spray manifold  $(M, {\bf G})$. The S-curvature ${\bf S}= {\bf S}_{(G, dV)} $ of $(\mathbf{G}, dV)$ is expressed by
\begin{equation}\label{621}
\mathbf{S}:=\frac{\partial G^{m}}{\partial y^{m}}-y^{m}\frac{\partial }{\partial x^{m}}(\ln \sigma).
\end{equation}

There is a mysterious
non-Riemannian quantities $\chi=\chi_{i}dx^{i}$. It is initially defined by the S-curvature ${\bf S}$  in \eqref{621},
but also given by the curvature $T^i_{\ j}$ in \eqref{Tij}
\begin{equation}\label{chi}
\chi_{i}=-\frac{1}{3}T^{m}_{~~i\cdot m}=\frac{1}{2}\{\mathbf{S}_{\cdot i|m}y^{m}-\mathbf{S}_{|i}\},
\end{equation}
where the covariant derivative $``|"$ is taken respect to $\mathbf{G}$.

Using ${\bf S}={\bf S}_{(G, dV)}$, one can modify the spray ${\bf G}$ to a new spray $\hat{\bf G}$ by
\begin{equation}
  \hat{\bf G}:= {\bf G}+\frac{2{\bf S}}{n+1} {\bf Y},\label{GGG}
\end{equation}
where ${\bf Y}= y^i\frac{\partial}{\partial y^{i}}$ is the canonical vertical vector field. $\hat{\bf G}$ is a globally defined spray on $M$. We call $\hat{\bf G}$ the projective spray associated with
$\bf G$ on $(M,dV)$. In a standard local coordinates in $TM$, $\hat{\bf G}=y^{i}\frac{\partial}{\partial x^{i}}-2 \hat{G}^{i}\frac{\partial}{\partial y^{i}}$ is given by
\[ \hat{G}^i = G^i -\frac{\bf S}{n+1} y^i.\]

By a direct computation, it is easy to see that  the S-curvature $\hat{\bf S}={\bf S}_{(\hat{\bf G}, dV)}$ always vanishes \cite{S1}. Thus the $\chi$-curvature $\hat{\chi}$ of $\hat{\bf G}$ vanishes by \eqref{chi}. That is, for any $(\mathbf G,dV)$,  the following holds
\begin{equation}\label{Shat}
\hat{\bf S}=0,~~ \hat{\chi}=0.
\end{equation}

Let $N^{i}_{k}$ and $\hat{N}^{i}_{k}$ be the non-linear connection coefficients of spray $\bf{G}$ and $\hat{\bf G}$  and
the Berwald connection coefficients are denoted by $\varGamma^{i}_{jk}$ and $\hat{\varGamma}^{i}_{jk}$, respectively.
It is easy to check that
\begin{eqnarray}
\hat{N}^{i}_{k}&=&N^{i}_{k}-\frac{\mathbf{S}}{n+1}\delta^{i}_{k}-\frac{\mathbf{S}_{\cdot k}}{n+1}y^{i} \label{a2}\\
\hat{\varGamma}^{i}_{jk}&=&\varGamma^{i}_{jk}-\frac{\mathbf{S}_{\cdot k}}{n+1}\delta^{i}_{j}
-\frac{\mathbf{S}_{\cdot j}}{n+1}\delta^{i}_{k}-\frac{\mathbf{S}_{\cdot k\cdot j}}{n+1}y^{i}.\nonumber
\end{eqnarray}
Then by \eqref{Hframe} and \eqref{a2}, we get
\begin{equation}\label{hat-delta}
 \hat{\frac{\delta}{\delta x^{k}}}=\frac{\delta}{\delta x^{k}}+\frac{\mathbf{S}_{\cdot k}}{n+1}\mathbf{Y}
 +\frac{\mathbf{S}}{n+1}\frac{\partial}{\partial y^{k}},
\end{equation}
where $\hat{\frac{\delta}{\delta x^{k}}}$ is defined in \eqref{Hframe} by using $\hat{\mathbf{G}}$.\\

For the Riemann curvatures of $\mathbf{G}$ and $\hat{\mathbf{G}}$, we have the following
\begin{lemma}\label{R1}(\cite{S4})
The Riemann curvatures and Ricci-scalar of  $\mathbf{G}$ and $\hat{\mathbf{G}}$ are related by
\begin{eqnarray}
\hat{R}^{i}_{~k}&=&R^{i}_{~k}+\tau\delta^{i}_{k}-\frac{1}{2}\tau_{\cdot k}+\frac{3\chi_{k}}{n+1}y^{i}\nonumber\\
\hat{R}&=&R+\tau,\label{Rtau}
\end{eqnarray}
where $\tau=(\frac{1}{n+1}\mathbf{S})^{2}+\frac{1}{n+1}\mathbf{S}_{|m}y^{m}$ and
$\chi_{k}=\frac{1}{2}(\mathbf{S}_{|m\cdot k}y^{m}-\mathbf{S}_{|k})$ and the covariant derivative $``|"$
 is taken respect to $\mathbf{G}$.
\end{lemma}

Let $f$ be a function defined on $TM$. From \eqref{hat-delta}, we have
\begin{equation}\label{hd}
f_{\parallel k}=f_{|k}+\frac{\mathbf{Y}(f)}{n+1}\mathbf{S}_{\cdot k}
 +\frac{\mathbf{S}}{n+1}f_{\cdot k},
\end{equation}
where the covariant derivative $``|"$  and $``\parallel"$ is taken respect to $\mathbf{G}$ and $\hat{\mathbf{G}}$ respectively.
Meanwhile, the Ricci identities  are given by
\begin{eqnarray}
&&f_{|m\cdot k}=f_{\cdot k|m}\label{Ri1},\\
&&f_{|k|m}-f_{|m|k}=f_{\cdot l}R^{l}_{~km},\label{Ri2}
\end{eqnarray}
where $R^{l}_{~km}$ is defined in \eqref{R3eq}.
Contracting \eqref{Ri2} with $y^{m}$, we obtain
\begin{equation}\label{Ri3}
f_{|k|0}-f_{|0|k}=f_{\cdot l}R^{l}_{~k}.
\end{equation}

The Weyl-type curvature of $\hat{\bf G}$ is given by
\[ \hat{T}^i_k = \hat{R}^i_{\ k} - (\hat{R}\delta^i_k -\frac{1}{2} \hat{R}_{\cdot k}y^i).\]
It turns out that $\hat{T}^i_k$ is independent of the volume form! Therefore it is a projective invariant. In fact, it is the well-known Weyl curvature
\begin{equation}\label{DefWik}
 W^{i}_{~k}=\hat{T}^{i}_{~k}.
\end{equation}
See Proposition 6.4 in \cite{S4}.
 Furthermore, one can rewrite \eqref{DefWik} as
\begin{eqnarray}
 \nonumber W^{i}_{~k}&=&T^{i}_{~k}+\frac{3\chi_{k}}{n+1}y^{i}\\
 &=&R^{i}_{~k}-(R\delta^{i}_{k}-\frac{1}{2}R_{\cdot k}y^{i})+\frac{3\chi_{k}}{n+1}y^{i}.\label{d7}
\end{eqnarray}
(\ref{d7}) is given in Lemma 3.3 in \cite{S4}.

The  Weyl-curvature $\mathbf{W}$ can be used to  characterize the sprays of scalar curvature.
\begin{proposition}(\cite{Z.Sza}\cite{S1})
Let $\bf{G}$ be a spray on a manifold. Then $\bf{G}$ is of scalar curvature if and only if  $\mathbf{W}=0$.
\end{proposition}

The following result proved by using the  Bianchi identity will be used later.
\begin{proposition}
For any spray $\mathbf{G}$ on a $n$-dimensional manifold $M$, we have the following identity:
\begin{equation}\label{New1}
R_{|p}-\frac{1}{2}R_{\cdot p|m}y^{m}-\frac{1}{n-1}(\chi_{p|m}y^{m}+R^{m}_{~p|m})=0,
\end{equation}
where the horizontal covariant derivative $``|"$ is taken respect to $\mathbf{G}$.
\end{proposition}
{\it Proof}:
The Bianchi identity says that
\begin{equation*}
(R^{i}_{jkl|p}+B^{i}_{jpm}R^{m}_{~kl})+(R^{i}_{~jlp|k}+B^{i}_{jlm}R^{m}_{~pk})+(R^{i}_{jpk|l}+B^{i}_{jkm}R^{m}_{~lp})=0.
\end{equation*}
Contracting the above equation with $y^{j}$ and using $y^{j}B^{i}_{jpm}=0,~y^{i}_{|p}=0$ and $R^{i}_{~kl}=y^{j}R^{i}_{jkl}$,
it follows that
\begin{equation*}
R^{i}_{~kl|p}+R^{i}_{~lp|k}+R^{i}_{~pk|l}=0.
\end{equation*}
Contracting it with $y^{l}$ and using $R^{i}_{~k}=R^{i}_{~kl}y^{l}$ yields that
\begin{equation}\label{b16}
R^{i}_{~k|p}-R^{i}_{~p|k}+R^{i}_{~pk|l}y^{l}=0.
\end{equation}
Then substituting \eqref{R3eq} into \eqref{b16}, it follows  that
\begin{equation}\label{Rikp}
R^{i}_{~k|p}-R^{i}_{~p|k}+\frac{1}{3}(R^{i}_{~p\cdot k}-R^{i}_{~k\cdot p})_{|l}y^{l}=0.
\end{equation}
Taking a trace of \eqref{Rikp} over $i$ and $k$, we get
\begin{equation}\label{Rikp1}
R_{|p}-\frac{1}{3}R_{\cdot p|m}y^{m}=\frac{1}{n-1}(R^{m}_{~p|m}-\frac{1}{3}R^{m}_{~p\cdot m}).
\end{equation}
Recall that the $\chi$-curvature of $\bf{G}$ is given by
\begin{equation*}
\chi_{i}=-\frac{1}{6}[2R^{m}_{~~i\cdot m}+(n-1)R_{\cdot{i}}]
\end{equation*}
which means
\begin{equation}\label{Rmjm}
R^{m}_{~~i\cdot m}=-3\chi_{i}-\frac{n-1}{2}R_{\cdot{i}}.
\end{equation}
Substituting \eqref{Rmjm} into equality \eqref{Rikp1}, we get the result.
\qed

\section{The Berwald-Weyl curvature $\mathbf{W}^{o}$}
In this section, we shall discuss the Berwald-Weyl curvature.
Let $(\mathbf{G}, dV)$ be a spray on an $n$-manifold $M$.  Let $\hat{\mathbf{G}}$ be the associated spray defined in (\ref{GGG}).
For a vector $y\in T_{x}M\setminus \{0\}$, the Berwald-Weyl curvature $\mathbf{W}^{o}_{y}:T_{x}M\rightarrow R$ is defined by
\[\mathbf{W}^{o}_{y}(u):=W^{o}_{k}(y)u^{k},~~~u=u^{k}\frac{\partial}{\partial x^{k}},\]where
\begin{equation}\label{W0k}
W^{o}_{k}:=\hat{R}_{\parallel k}-\frac{1}{2}\hat{R}_{\cdot k\parallel m}y^{m}.
\end{equation}
Here $\hat{R}=\frac{1}{n-1}\hat{R}^{m}_{~m}$ denotes the Ricci scalar of $\hat{\mathbf{G}}$ and the covariant derivative
$``\parallel"$ is taken respect to $\hat{\mathbf{G}}$.
By Ricci identity \eqref{Ri1}, we can rewrite $W^o_{k}$ as
\begin{equation}\label{W0k1}
W^{o}_{k}=\frac{1}{2}\{3\hat{R}_{\parallel k}-(\hat{R}_{\parallel m}y^{m})_{\cdot k}\}.
\end{equation}

Since $\hat{\mathbf{G}}$ is a projective invariant with respect to a fixed volume form $dV$,  $\mathbf{W}^{o}$ is a projective invariant.
By the homogeneity of $\hat{R}$ and $y^{k}_{\parallel m}=0$, it follows from \eqref{W0k} or \eqref{W0k1} that
\begin{equation*}
W^{o}_{k}y^{k}=0.
\end{equation*}

Applying formula \eqref{New1} to $\hat{\bf{G}}$ and using \eqref{Shat}, we obtain
\begin{equation}\label{Bweq1}
\hat{R}_{||k}-\frac{1}{2}\hat{R}_{\cdot p||m}y^{m}=\frac{1}{n-1}\hat{R}^{m}_{~k||m},
\end{equation}
where the covariant derivative $``\parallel"$ is taken respect to $\hat{\mathbf{G}}$ and $\hat{R}^{m}_{~k}$ is the Riemann curvature of $\hat{\mathbf{G}}$. By \eqref{W0k} and \eqref{Bweq1}, we immediately obtain another expression for  $W^{o}_{k}$
\[W^{o}_{k}=\frac{1}{n-1}\hat{R}^{m}_{~k||m}.\]

\vspace{0.5cm}
\noindent
{\it Proof of Proposition \ref{BWC1}}:
Firstly, according to \eqref{Rtau} and \eqref{hd}, one gets
\begin{eqnarray}
\hat{R}_{\parallel k}&=&(R+\tau)_{|k}+\frac{2(R+\tau)}{n+1}\mathbf{S}_{\cdot k}+\frac{(R+\tau)_{\cdot k}}{n+1}\mathbf{S}\label{R}\\
\hat{R}_{\parallel m} y^{m}&=&(R+\tau)_{|m}y^{m}+\frac{4(R+\tau)}{n+1}\mathbf{S}.\label{Rk2}
\end{eqnarray}
Taking the derivative of \eqref{Rk2} with respect to $y^{k}$ yields that
\begin{eqnarray}\label{Rky}
(\hat{R}_{\parallel m} y^{m})_{\cdot k}=(R+\tau)_{|k}+(R+\tau)_{|m\cdot k}y^{m}
+\frac{4(R+\tau)}{n+1}\mathbf{S}_{\cdot k}+\frac{4(R+\tau)_{\cdot k}}{n+1}\mathbf{S}.
\end{eqnarray}
Plugging \eqref{R} and \eqref{Rky} into \eqref{W0k1}, we obtain
\begin{eqnarray}\label{W0k2}
W^{o}_{k}&=&R_{|k}-\frac{1}{2}R_{\cdot k|m}y^{m}-\frac{1}{2(n+1)}[\mathbf{S}_{|l\cdot k|m}y^{m}y^{l}-\mathbf{S}_{|k|m}y^{m}]\nonumber\\
&-&\frac{3\mathbf{S}}{2(n+1)^{2}}(\mathbf{S}_{|m\cdot k}y^{m}-\mathbf{S}_{|k})-\frac{1}{n+1}[\mathbf{S}_{|k|m}y^{m}-\mathbf{S}_{|m|k}y^{m}]\nonumber\\
&+&\frac{R\mathbf{S}_{\cdot k}}{n+1}-\frac{\mathbf{S}R_{\cdot k}}{2(n+1)}.
\end{eqnarray}
Now, by using \eqref{chi} and \eqref{Ri2}, $W^{o}_{k}$ in \eqref{W0k2} reduces to
\begin{eqnarray}
\nonumber W^{o}_{k}&=&R_{|k}-\frac{1}{2}R_{\cdot k|m}y^{m}-\frac{1}{n+1}\chi_{k|m}y^{m}-\frac{3\mathbf{S}}{(n+1)^{2}}\chi_{k}\\
&-&\frac{1}{n+1}R^{l}_{k}\mathbf{S}_{\cdot l}+\frac{R\mathbf{S}_{\cdot k}}{n+1}-\frac{\mathbf{S}R_{\cdot k}}{2(n+1)}\label{d5}
\end{eqnarray}
The last three terms in \eqref{d5} are rewritten by
\begin{eqnarray}
&&-\frac{1}{n+1}R^{l}_{k}\mathbf{S}_{\cdot l}+\frac{R\mathbf{S}_{\cdot k}}{n+1}-\frac{\mathbf{S}R_{\cdot k}}{2(n+1)}\nonumber\\
&=&-\frac{1}{n+1}[R^{l}_{k}-(R\delta^{l}_{k}-\frac{1}{2}R_{\cdot k}y^{l})]\mathbf{S}_{\cdot l}\nonumber\\
&=&-\frac{1}{n+1}T^{l}_{k}\mathbf{S}_{\cdot l},\label{d6}
\end{eqnarray}
where $T^{l}_{k}$ is defined in \eqref{Tij}. Taking \eqref{d6} into \eqref{d5}, we obtain
\begin{equation}\label{b2}
W^{o}_{k}=R_{|k}-\frac{1}{2}R_{\cdot k|m}y^{m}-\frac{1}{n+1}\chi_{k|j}y^{j}-\frac{3\mathbf{S}}{(n+1)^{2}}\chi_{k}
-\frac{1}{n+1}T^{l}_{k}\mathbf{S}_{\cdot l}.
\end{equation}
Finally, submitting \eqref{d7} into \eqref{b2}, we obtain
\begin{equation}\label{d99**}
W^{o}_{k}=R_{|k}-\frac{1}{2}R_{\cdot k|m}y^{m}-\frac{1}{n+1}\chi_{k|m}y^{m}-\frac{1}{n+1}W^{m}_{~k}\mathbf{S}_{\cdot m},
\end{equation}
This completes the proof.
\qed \\

Note that if a Finsler metric $F$ has constant S-curvature, that is, the S-curvature of F with respect to $dV=dV_{BH}$ satisfies  that ${\bf S}= (n+1) c F$, $c=constant$, then $\chi=0$.
By \eqref{d99**}, we obtain  Theorem \ref{E-RC}.
We  can generalize Theorem \ref{E-RC} as follows.
\begin{proposition}\label{ConS-Ein}
Let $F$ be a Finsler metric on an $n$-dimensional manifold $M$. If $F$ is of constant S-curvature and is Einstein with $\mathbf{Ric}=(n-1)\sigma F^{2}$, then
\[W^{o}_{k}=F^{3}(\frac{\theta}{F})_{\cdot k},\]
for the $dV=dV_{BH}$, where $\theta=\sigma_{x^{m}}y^{m}.$
\end{proposition}

\noindent
{\it Proof}:
By assumption, we have
\begin{equation}\label{Ein}
\mathbf{S}=(n+1)cF,~~ R=\sigma F^{2},
\end{equation}
where $c$ is a constant and $\sigma=\sigma(x)$ is a function on $M$.
Then by \eqref{chi}, we have
\begin{equation}\label{chiR}
\chi=0,~~\mathbf{S}_{\cdot m}=(n+1)cF^{-1}y_{m},
\end{equation}
where $y_{m}:=\frac{1}{2}(F^{2})_{\cdot m}$.
By using \eqref{chiR} and \eqref{d7}, we get
\begin{equation}\label{term4}
W^{m}_{~k}\mathbf{S}_{\cdot m}=-(R\mathbf{S}_{\cdot k}-\frac{1}{2}\mathbf{S}R_{\cdot k})=0,
\end{equation}
where we use $R^{m}_{~k}y_{m}$=0.
On the other hand, by a direct calculation, we obtain
\begin{equation}\label{term12}
R_{|k}-\frac{1}{2}R_{\cdot k|m}y^{m}=F^{3}(\frac{\theta}{F})_{\cdot k}
\end{equation}
with $\theta=\sigma_{x^{k}}y^{k}$. Plugging \eqref{term4} and \eqref{term12} into \eqref{d99**}, we complete the proof.
\qed

\vspace{0.5cm}

It is known that if a Finsler metric is of constant flag curvature, then $\chi=0$. By Corollary \ref{cor1.4}, we have
\begin{cor}
For any two-dimensional Finsler metric of constant flag curvature, ${\bf W}^o=0$.
\end{cor}

\vspace{0.5cm}
\noindent
{\it Proof of Proposition \ref{W0kN3}}:
Taking the horizontal covariant derivatives $``|"$  with respect to $\mathbf{G}$ of \eqref{d7}, we have
\begin{equation}\label{Whori}
W^{m}_{k|m}=R^{m}_{k|m}-(R_{|k}-\frac{1}{2}R_{\cdot k|m}y^{m})+\frac{3}{n+1}\chi_{k|m}y^{m}.
\end{equation}
By \eqref{New1}, it follows that
\begin{equation}\label{New11}
R^{m}_{k|m}=(n-1)(R_{|k}-\frac{1}{2}R_{\cdot k|m}y^{m})-\chi_{k|m}y^{m}.
\end{equation}
Substituting \eqref{New11} into \eqref{Whori} yields that
\begin{equation}\label{Wmkm}
W^{m}_{k|m}=(n-2)(R_{|k}-\frac{1}{2}R_{\cdot k|m}y^{m}-\frac{1}{n+1}\chi_{k|m}y^{m}).
\end{equation}
By \eqref{Wmkm} and the assumption $n\geq3$, $W^{0}_{k}$ in \eqref{d99**} becomes
\begin{eqnarray}
W^{0}_{k}&=&\frac{1}{n-2}W^{m}_{k|m}-\frac{1}{n+1}W^{m}_{k}\mathbf{S}_{\cdot m}\label{d91}\\
&=&\frac{1}{n-2}W^{m}_{~k||m}\nonumber,
\end{eqnarray}
where the covariant derivative $``\parallel"$ is taken respect to $\hat{\mathbf{G}}$.
\qed

\vspace{0.5cm}
We have another method to prove Proposition \ref{W0kN3}. Firstly, taking the horizontal covariant derivatives $``|"$  with respect to $\mathbf{G}$ from \eqref{Tij}, we obtain
\begin{equation*}
T^{m}_{~k|m}=R^{m}_{k|m}-(R_{|k}-\frac{1}{2}R_{\cdot k|m}y^{m}).
\end{equation*}
Plugging \eqref{New11} into preceding equation yields that
\begin{equation}\label{Tmjm}
T^{m}_{k|m}=(n-2)(R_{|k}-\frac{1}{2}R_{\cdot k|m}y^{m})-\chi_{k|m}y^{m}.
\end{equation}
Let $\hat{T}^{i}_{k}$ be the Weyl-type curvature of projective spray $\hat{\mathbf{G}}$. Applying \eqref{Tmjm} to $\hat{\mathbf{G}}$ and using \eqref{Shat}, we obtain
\begin{equation}\label{hatTmk}
\hat{T}^{m}_{k\parallel m}=(n-2)(\hat{R}_{\parallel k}-\frac{1}{2}\hat{R}_{\cdot k\parallel m}y^{m}),
\end{equation}
where the covariant derivative $``\parallel"$ is taken respect to $\hat{\mathbf{G}}$.
Therefore, according to \eqref{W0k}, \eqref{hatTmk} and \eqref{DefWik}, we obtain
\begin{equation*}
W^{o}_{k}=\hat{R}_{\parallel k}-\frac{1}{2}\hat{R}_{\cdot k\parallel m}y^{m}=\frac{1}{n-2}\hat{T}^{m}_{k\parallel m}=\frac{1}{n-2}W^{m}_{k||m}.
\end{equation*}

\section{BWeyl-flat}
In this section, we shall define the notion of BWeyl-flat sprays and give three equivalent conditions for a spray to be BWeyl-flat when $\dim M\geq3$.

Let $\mathbf{G}$ be a spray on an $n$-dimensional manifold $M$ of $n\geq3$. By \eqref{d99**}, we know that the Berwald-Weyl curvature $\mathbf{W}^{o}$ of $\mathbf{G}$ depends on the volume form $dV$ in general.
We make the following
\begin{definition}\label{WeylF}A spray  ${\bf G}$ on an $n$-dimensional manifold $M$ of $n\geq3$ is said to be BWeyl-flat if there is a volume
form $dV$ on $M$ such that Berwald-Weyl curvature $\mathbf{W}^{o}$ vanishes, that is
\[ \bf{W}^{o}=0.\]
 A Finsler metric $F$ on $M$ is said to be BWeyl-flat if the induced spray ${\bf G}={\bf G}_F$ is BWeyl-flat.
\end{definition}
\begin{remark}
 We have obtained from Theorem\ref{e88} that a spray of scalar curvature must have $\bf{W}^{o}=0$ for any volume form. Therefore, a spray of scalar curvature must be BWeyl-flat in our sense.
\end{remark}

Let $dV$ and $d\tilde{V}$ be volume forms on $M$ with $dV = e^{(n+1)f} d\tilde{V}$, where $f=f(x)$ is a scalar function on $M$.
Then the S-curvatures ${\bf S}= {\bf S}_{(G, dV)}$ and $\tilde{\bf S}= {\bf S}_{(G, d\tilde{V})}$ are related by
\begin{equation}\label{516}
\mathbf{S}
=\tilde{\mathbf{S}}-(n+1) f_0,
\end{equation}
where $f_0 := f_{x^m}(x)y^m$.
By (\ref{516}), one sees that ${\bf S}_{(G, dV)} =(n+1) \phi_0$ if and only if ${\bf S}_{(G, d\tilde{V})} =(n+1) \tilde{\phi}_0$
with $\phi_0 = \tilde{\phi}_0 -f_0$.

Let $W^{o}_{k}$ and $\widetilde{W}^{o}_{k}$ be the Berwald-Weyl curvature of $\mathbf{G}$ with respect to volume
form $dV$ and $d\tilde{V}$ respectively. It follows from \eqref{d99**} that
\begin{eqnarray}
W^{o}_{k}&=&R_{|k}-\frac{1}{2}R_{\cdot k|m}y^{m}-\frac{1}{n+1}\chi_{k|m}y^{m}-\frac{1}{n+1}W^{m}_{k}\mathbf{S}_{\cdot m}\label{eqP1},\\
\widetilde{W}^{o}_{k}&=&R_{|k}-\frac{1}{2}R_{\cdot k|m}y^{m}-\frac{1}{n+1}\chi_{k|m}y^{m}-\frac{1}{n+1}W^{m}_{k}\mathbf{\tilde{S}}_{\cdot m},\label{eqP2}
\end{eqnarray}
where $``|"$ is the horizontal covariant derivative with respect to $\mathbf{G}$.
It follows from \eqref{516}, \eqref{eqP1} and \eqref{eqP2} that
\begin{equation}\label{5297}
\widetilde{W}^{o}_{k}=W^{o}_{k}-W^{m}_{k}f_{m}.
\end{equation}

Now we give some equivalent conditions for a spray to be BWeyl-flat
in $dimM\geq3$.
\begin{theorem}\label{5310}
Let  ${\bf G}$ be a spray  on an $n$-dimensional manifold $M$ of $n\geq 3$. Then the following are equivalent:
\ben
\item[(a)]  $\mathbf{G}$ is of BWeyl-flat.
\item[(b)]  for any  volume form $dV$ on $M$ there is a scalar function $f$ on $M$ such that
\be
W^{o}_{k}=W^{m}_{k}f_{m}.\label{oPRicf}
\ee
\item[(c)]  for any  volume form $dV$ on $M$ there is  a scalar function $f$ on $M$ such that
\be
W^{m}_{~k|m}=(n-2)W^{m}_{k}\Xi_{\cdot m},\label{Wmkm1}
\ee
\een
where
$``|"$ is the horizontal covariant derivative with respect to $\mathbf{G}$ and $\Xi := \frac{{\bf S}}{n+1}+f_0 $
with $f_{0}=f_{x^{m}}y^{m}$ and ${\bf S}= {\bf S}_{(G, dV)}$.
\end{theorem}
{\it Proof}:
(a) $\Rightarrow$ (b).   Assume that for some volume form $d\tilde{V}$,
\[\mathbf{\widetilde{W}}^{o}=0.\]
Then by (\ref{5297}), for any volume form  $dV= e^{(n+1)f} d\tilde{V}$,
\begin{equation}
W^{o}_{k}=W^{m}_{k}f_{m}. \label{PRic*}
\end{equation}

(b) $\Rightarrow$ (c).
By assumption and (\ref{d91}), we have
\begin{eqnarray}
W^{m}_{~k|m}&=&\frac{n-2}{n+1}W^{m}_{~k} \mathbf{S}_{\cdot m}+(n-2)W^{m}_{~k}f_{m}\label{atob}\\
&=&(n-2)W^{m}_{~k}\Xi_{\cdot m},\nonumber
\end{eqnarray}
where $\Xi = \frac{\bf S}{n+1}+f_0$.

(c) $\Rightarrow$ (a).
Let $dV$ be an arbitrary  volume form on $M$. There is a scalar function $f $ on $M$ such that (\ref{Wmkm1}) holds. That is
\[ W^{m}_{~k|m}=(n-2)W^{m}_{k}\Xi_{\cdot m}.\]
Taking
$d\tilde{V}= e^{-(n+1) f} dV$,  we have
\[ \tilde{\bf S} = {\bf S} +(n+1) f_0= (n+1) \Xi.\]
By (\ref{d91}), (\ref{5297}) and \eqref{atob}, we obtain
\[ \widetilde{W}^{o}_{k}=W^{o}_{k}-W^{m}_{k}f_{m}
=\frac{1}{n-2}W^{m}_{k|m}-\frac{1}{n+1}W^{m}_{k}\mathbf{S}_{\cdot m}-W^{m}_{k}f_{m}=0.\]

The condition on the S-curvature ${\bf S}_{(G, dV)}=(n+1) \phi_0 $ being an exact 1-form, is actually a condition on the spray, independent of the choice of a particular volume form $dV$. If a spray ${\bf G}$ satisfies that ${\bf S}_{(G, dV)}=(n+1) \phi_0 $ for some scalar function $\phi$ on $M$, then $\Xi=0$ for $f=-\phi$ holds, where $\Xi = \frac{\bf S}{n+1}+f_0$. Thus in view of \eqref {Wmkm1}, we obtain the following

\begin{cor}
Let ${\bf G}$ be a spray on an $n$-dimensional manifold $M$ of $n\geq 3$.
Suppose that $W^{m}_{~k|m}=0$ and ${\bf S}_{(G, dV)}=(n+1) \phi_0 $ for some scalar function $\phi$, then ${\bf G}$ is BWeyl-flat.
\end{cor}

A spray is said to be projectively Ricci-flat if there is a volume form $dV$ on $M$ such that $\hat{R}=\frac{1}{n-1}\mathbf{Ric}_{\hat{G}}=0$.
A Finsler metric $F$ on $M$ is said to be projectively Ricci-flat if the induced spray $\mathbf{G}=\mathbf{G}_{F}$ is projectively Ricci-flat\cite{S5}. Thus from \eqref{7261}, if a spray or a Finsler metric is projectively Ricci-flat, then it is BWeyl-flat.

\begin{example}
Let $\alpha:=(1+4|x|^{2})\sqrt{\frac{(1+4|x|^{2})|y|^2-4\langle x,y\rangle}{1+4|x|^{2}}}$ and $\beta=2\langle x,y\rangle$. It is shown that the square metric $F=\frac{(\alpha+\beta)^{2}}{\alpha}$ is projectively Ricci-flat\cite{Zhu}. Thus $F$ is BWeyl-flat. On the other hand, by using Theorem 1.2 in \cite{CS3}, one can see that $F$ is not of isotropic S-curvature by a direct computation.
\end{example}

\vspace{0.6cm}

\noindent Zhongmin Shen

\noindent Department of Mathematical Sciences,

\noindent  Indiana University-Purdue University Indianapolis, IN 46202-3216, USA.

\noindent \verb"zshen@iupui.edu"

\vspace{0.6cm}

\noindent Liling Sun

\noindent School of Science,

\noindent Jimei University, Xiamen 361021, P.R. China.

\noindent \verb"sunliling@jmu.edu.cn"

\end{document}